

\input amstex
\expandafter\ifx\csname mathdefs.tex\endcsname\relax
  \expandafter\gdef\csname mathdefs.tex\endcsname{}
\else \message{Hey!  Apparently you were trying to
  \string\input{mathdefs.tex} twice.   This does not make sense.} 
\errmessage{Please edit your file (probably \jobname.tex) and remove
any duplicate ``\string\input'' lines}\endinput\fi




\catcode`\X=12\catcode`\@=11

\def\n@wcount{\alloc@0\count\countdef\insc@unt}
\def\n@wwrite{\alloc@7\write\chardef\sixt@@n}
\def\n@wread{\alloc@6\read\chardef\sixt@@n}
\def\r@s@t{\relax}\def\v@idline{\par}\def\@mputate#1/{#1}
\def\l@c@l#1X{\firstpart.#1}\def\gl@b@l#1X{#1}\def\t@d@l#1X{{}}

\def\crossrefs#1{\ifx\all#1\let\tr@ce=\all\else\def\tr@ce{#1,}\fi
   \n@wwrite\cit@tionsout\openout\cit@tionsout=\jobname.cit 
   \write\cit@tionsout{\tr@ce}\expandafter\setfl@gs\tr@ce,}
\def\setfl@gs#1,{\def\@{#1}\ifx\@\empty\let\next=\relax
   \else\let\next=\setfl@gs\expandafter\xdef
   \csname#1tr@cetrue\endcsname{}\fi\next}
\def\m@ketag#1#2{\expandafter\n@wcount\csname#2tagno\endcsname
     \csname#2tagno\endcsname=0\let\tail=\all\xdef\all{\tail#2,}
   \ifx#1\l@c@l\let\tail=\r@s@t\xdef\r@s@t{\csname#2tagno\endcsname=0\tail}\fi
   \expandafter\gdef\csname#2cite\endcsname##1{\expandafter
     \ifx\csname#2tag##1\endcsname\relax?\else\csname#2tag##1\endcsname\fi
     \expandafter\ifx\csname#2tr@cetrue\endcsname\relax\else
     \write\cit@tionsout{#2tag ##1 cited on page \folio.}\fi}
   \expandafter\gdef\csname#2page\endcsname##1{\expandafter
     \ifx\csname#2page##1\endcsname\relax?\else\csname#2page##1\endcsname\fi
     \expandafter\ifx\csname#2tr@cetrue\endcsname\relax\else
     \write\cit@tionsout{#2tag ##1 cited on page \folio.}\fi}
   \expandafter\gdef\csname#2tag\endcsname##1{\expandafter
      \ifx\csname#2check##1\endcsname\relax
      \expandafter\xdef\csname#2check##1\endcsname{}%
      \else\immediate\write16{Warning: #2tag ##1 used more than once.}\fi
      \multit@g{#1}{#2}##1/X%
      \write\t@gsout{#2tag ##1 assigned number \csname#2tag##1\endcsname\space
      on page \number\count0.}%
   \csname#2tag##1\endcsname}}
\def\multit@g#1#2#3/#4X{\def\t@mp{#4}\ifx\t@mp\empty%
      \global\advance\csname#2tagno\endcsname by 1 
      \expandafter\xdef\csname#2tag#3\endcsname
      {#1\number\csname#2tagno\endcsnameX}%
   \else\expandafter\ifx\csname#2last#3\endcsname\relax
      \expandafter\n@wcount\csname#2last#3\endcsname
      \global\advance\csname#2tagno\endcsname by 1 
      \expandafter\xdef\csname#2tag#3\endcsname
      {#1\number\csname#2tagno\endcsnameX}
      \write\t@gsout{#2tag #3 assigned number \csname#2tag#3\endcsname\space
      on page \number\count0.}\fi
   \global\advance\csname#2last#3\endcsname by 1
   \def\t@mp{\expandafter\xdef\csname#2tag#3/}%
   \expandafter\t@mp\@mputate#4\endcsname
   {\csname#2tag#3\endcsname\lastpart{\csname#2last#3\endcsname}}\fi}
\def\t@gs#1{\def\all{}\m@ketag#1e\m@ketag#1s\m@ketag\t@d@l p
   \m@ketag\gl@b@l r \n@wread\t@gsin
   \openin\t@gsin=\jobname.tgs \re@der \closein\t@gsin
   \n@wwrite\t@gsout\openout\t@gsout=\jobname.tgs }
\outer\def\localtags{\t@gs\l@c@l}
\outer\def\globaltags{\t@gs\gl@b@l}
\outer\def\newlocaltag#1{\m@ketag\l@c@l{#1}}
\outer\def\newglobaltag#1{\m@ketag\gl@b@l{#1}}

\newif\ifpr@ 
\def\m@kecs #1tag #2 assigned number #3 on page #4.%
   {\expandafter\gdef\csname#1tag#2\endcsname{#3}
   \expandafter\gdef\csname#1page#2\endcsname{#4}
   \ifpr@\expandafter\xdef\csname#1check#2\endcsname{}\fi}
\def\re@der{\ifeof\t@gsin\let\next=\relax\else
   \read\t@gsin to\t@gline\ifx\t@gline\v@idline\else
   \expandafter\m@kecs \t@gline\fi\let \next=\re@der\fi\next}
\def\pretags#1{\pr@true\pret@gs#1,,}
\def\pret@gs#1,{\def\@{#1}\ifx\@\empty\let\n@xtfile=\relax
   \else\let\n@xtfile=\pret@gs \openin\t@gsin=#1.tgs \message{#1} \re@der 
   \closein\t@gsin\fi \n@xtfile}

\newcount\sectno\sectno=0\newcount\subsectno\subsectno=0
\newif\ifultr@local \def\ultralocal{\ultr@localtrue}
\def\firstpart{\number\sectno}
\def\lastpart#1{\ifcase#1 \or a\or b\or c\or d\or e\or f\or g\or h\or 
   i\or k\or l\or m\or n\or o\or p\or q\or r\or s\or t\or u\or v\or w\or 
   x\or y\or z \fi}

\def\resetall{\global\advance\sectno by 1\subsectno=0
   \gdef\firstpart{\number\sectno}\r@s@t}
\def\resetsub{\global\advance\subsectno by 1
   \gdef\firstpart{\number\sectno.\number\subsectno}\r@s@t}
\def\newsection#1\par{\resetall\vskip0pt plus.3\vsize\penalty-250
   \vskip0pt plus-.3\vsize\bigskip\bigskip
   \message{#1}\leftline{\bf#1}\nobreak\bigskip}
\def\subsection#1\par{\ifultr@local\resetsub\fi
   \vskip0pt plus.2\vsize\penalty-250\vskip0pt plus-.2\vsize
   \bigskip\smallskip\message{#1}\leftline{\bf#1}\nobreak\medskip}

\def\t@gsoff#1,{\def\@{#1}\ifx\@\empty\let\next=\relax\else\let\next=\t@gsoff
   \def\@@{p}\ifx\@\@@\else
   \expandafter\gdef\csname#1cite\endcsname##1{\zeigen{##1}}
   \expandafter\gdef\csname#1page\endcsname##1{?}
   \expandafter\gdef\csname#1tag\endcsname##1{\zeigen{##1}}\fi\fi\next}
\def\verbatimtags{\ifx\all\relax\else\expandafter\t@gsoff\all,\fi}
\def\zeigen#1{\hbox{$\langle$}#1\hbox{$\rangle$}}

\def\(#1){\edef\dot@g{\ifmmode\ifinner(\hbox{\noexpand\etag{#1}})
   \else\noexpand\eqno(\hbox{\noexpand\etag{#1}})\fi
   \else(\noexpand\ecite{#1})\fi}\dot@g}

\newif\ifbr@ck
\def\eat#1{}
\def\[#1]{\br@cktrue[\br@cket#1'X]}
\def\br@cket#1'#2X{\def\temp{#2}\ifx\temp\empty\let\next\eat
   \else\let\next\br@cket\fi
   \ifbr@ck\br@ckfalse\br@ck@t#1,X\else\br@cktrue#1\fi\next#2X}
\def\br@ck@t#1,#2X{\def\temp{#2}\ifx\temp\empty\let\neext\eat
   \else\let\neext\br@ck@t\def\temp{,}\fi
   \def\teemp{#1}\ifx\teemp\empty\else\rcite{#1}\fi\temp\neext#2X}
\def\resetbr@cket{\gdef\[##1]{[\rtag{##1}]}}
\def\references{\resetbr@cket\newsection References\par}

\newtoks\symb@ls\newtoks\s@mb@ls\newtoks\p@gelist\n@wcount\ftn@mber
    \ftn@mber=1\newif\ifftn@mbers\ftn@mbersfalse\newif\ifbyp@ge\byp@gefalse
\def\defm@rk{\ifftn@mbers\n@mberm@rk\else\symb@lm@rk\fi}
\def\n@mberm@rk{\xdef\m@rk{{\the\ftn@mber}}%
    \global\advance\ftn@mber by 1 }
\def\rot@te#1{\let\temp=#1\global#1=\expandafter\r@t@te\the\temp,X}
\def\r@t@te#1,#2X{{#2#1}\xdef\m@rk{{#1}}}
\def\b@@st#1{{$^{#1}$}}\def\str@p#1{#1}
\def\symb@lm@rk{\ifbyp@ge\rot@te\p@gelist\ifnum\expandafter\str@p\m@rk=1 
    \s@mb@ls=\symb@ls\fi\write\f@nsout{\number\count0}\fi \rot@te\s@mb@ls}
\def\byp@ge{\byp@getrue\n@wwrite\f@nsin\openin\f@nsin=\jobname.fns 
    \n@wcount\currentp@ge\currentp@ge=0\p@gelist={0}
    \re@dfns\closein\f@nsin\rot@te\p@gelist
    \n@wread\f@nsout\openout\f@nsout=\jobname.fns }
\def\m@kelist#1X#2{{#1,#2}}
\def\re@dfns{\ifeof\f@nsin\let\next=\relax\else\read\f@nsin to \f@nline
    \ifx\f@nline\v@idline\else\let\t@mplist=\p@gelist
    \ifnum\currentp@ge=\f@nline
    \global\p@gelist=\expandafter\m@kelist\the\t@mplistX0
    \else\currentp@ge=\f@nline
    \global\p@gelist=\expandafter\m@kelist\the\t@mplistX1\fi\fi
    \let\next=\re@dfns\fi\next}
\def\symbols#1{\symb@ls={#1}\s@mb@ls=\symb@ls} 
\def\bigsymbol{\textstyle}
\symbols{\bigsymbol\ast,\dagger,\ddagger,\sharp,\flat,\natural,\star}
\def\ftnumbers{\ftn@mberstrue} \def\ftsymbols{\ftn@mbersfalse}
\def\paginal{\byp@ge} \def\resetftnumbers{\ftn@mber=1}
\def\ftnote#1{\defm@rk\expandafter\expandafter\expandafter\footnote
    \expandafter\b@@st\m@rk{#1}}

\long\def\jump#1\endjump{}
\def\ssum{\mathop{\lower .1em\hbox{$\textstyle\Sigma$}}\nolimits}

\def\qed{\nobreak\kern 1em \vrule height .5em width .5em depth 0em}
\def\newneq{\hbox{\rlap{\hbox to 1\wd9{\hss$=$\hss}}\raise .1em 
   \hbox to 1\wd9{\hss$\scriptscriptstyle/$\hss}}}
\def\subsetne{\setbox9 = \hbox{$\subset$}\mathrel{\hbox{\rlap
   {\lower .4em \newneq}\raise .13em \hbox{$\subset$}}}}
\def\supsetne{\setbox9 = \hbox{$\subset$}\mathrel{\hbox{\rlap
   {\lower .4em \newneq}\raise .13em \hbox{$\supset$}}}}

\def\vbar{\mathchoice{\vrule height6.3ptdepth-.5ptwidth.8pt\kern-.8pt}
   {\vrule height6.3ptdepth-.5ptwidth.8pt\kern-.8pt}
   {\vrule height4.1ptdepth-.35ptwidth.6pt\kern-.6pt}
   {\vrule height3.1ptdepth-.25ptwidth.5pt\kern-.5pt}}
\def\f@dge{\mathchoice{}{}{\mkern.5mu}{\mkern.8mu}}
\def\b@c#1#2{{\rm \mkern#2mu\vbar\mkern-#2mu#1}}
\def\b@b#1{{\rm I\mkern-3.5mu #1}}
\def\b@a#1#2{{\rm #1\mkern-#2mu\f@dge #1}}
\def\bb#1{{\count4=`#1 \advance\count4by-64 \ifcase\count4\or\b@a A{11.5}\or
   \b@b B\or\b@c C{5}\or\b@b D\or\b@b E\or\b@b F \or\b@c G{5}\or\b@b H\or
   \b@b I\or\b@c J{3}\or\b@b K\or\b@b L \or\b@b M\or\b@b N\or\b@c O{5} \or
   \b@b P\or\b@c Q{5}\or\b@b R\or\b@a S{8}\or\b@a T{10.5}\or\b@c U{5}\or
   \b@a V{12}\or\b@a W{16.5}\or\b@a X{11}\or\b@a Y{11.7}\or\b@a Z{7.5}\fi}}

\catcode`\X=11 \catcode`\@=12


\expandafter\ifx\csname citeadd.tex\endcsname\relax
\expandafter\gdef\csname citeadd.tex\endcsname{}
\else \message{Hey!  Apparently you were trying to
\string\input{citeadd.tex} twice.   This does not make sense.} 
\errmessage{Please edit your file (probably \jobname.tex) and remove
any duplicate ``\string\input'' lines}\endinput\fi

\localtags
\NoBlackBoxes
\define\mr{\medskip\roster}
\define\sn{\smallskip\noindent}
\define\mn{\medskip\noindent}
\define\bn{\bigskip\noindent}
\define\ub{\underbar}

\define\ermn{\endroster\medskip\noindent}

\define \nl{\newline}
\magnification=\magstep 1
\documentstyle {amsppt}
\topmatter
\title{E20 - A continuation of [DjSh:691]} \endtitle
\author {Saharon Shelah \thanks {\null\newline I would like to thank 
Alice Leonhardt for the beautiful typing. \null\newline
 First Typed - 98/Nov/23 \null\newline
 Latest Revision - 99/Jan/19} \endthanks} \endauthor  
\affil{Institute of Mathematics\\
 The Hebrew University\\
 Jerusalem, Israel
 \medskip
 Rutgers University\\
 Mathematics Department\\
 New Brunswick, NJ  USA} \endaffil
\endtopmatter
\document  

\expandafter\ifx\csname alice2jlem.tex\endcsname\relax
  \expandafter\gdef\csname alice2jlem.tex\endcsname{}
\else \message{Hey!  Apparently you were trying to
\string\input{alice2jlem.tex}  twice.   This does not make sense.}
\errmessage{Please edit your file (probably \jobname.tex) and remove
any duplicate ``\string\input'' lines}\endinput\fi

\expandafter\ifx\csname bib4plain.tex\endcsname\relax
  \expandafter\gdef\csname bib4plain.tex\endcsname{}
\else \message{Hey!  Apparently you were trying to \string\input
  bib4plain.tex twice.   This does not make sense.}
\errmessage{Please edit your file (probably \jobname.tex) and remove
any duplicate ``\string\input'' lines}\endinput\fi

\def\renewcommand{\newcommand}	       
\edef\cite{\the\catcode`@}%
\catcode`@ = 11
\let\@oldatcatcode = \cite
\chardef\@letter = 11
\chardef\@other = 12
%
%
%
%
\def\@innerdef#1#2{\edef#1{\expandafter\noexpand\csname #2\endcsname}}%
%
%
\@innerdef\@innernewcount{newcount}%
\@innerdef\@innernewdimen{newdimen}%
\@innerdef\@innernewif{newif}%
\@innerdef\@innernewwrite{newwrite}%
%
%
%
\def\@gobble#1{}%
%
%
%
\ifx\inputlineno\@undefined
   \let\@linenumber = \empty 
\else
   \def\@linenumber{\the\inputlineno:\space}%
\fi
%
%
%
\def\@futurenonspacelet#1{\def\cs{#1}%
   \afterassignment\@stepone\let\@nexttoken=
}%
\begingroup 
\def\\{\global\let\@stoken= }%
\\ 
\endgroup
\def\@stepone{\expandafter\futurelet\cs\@steptwo}%
\def\@steptwo{\expandafter\ifx\cs\@stoken\let\@@next=\@stepthree
   \else\let\@@next=\@nexttoken\fi \@@next}%
\def\@stepthree{\afterassignment\@stepone\let\@@next= }%
%
%
%
\def\@getoptionalarg#1{%
   \let\@optionaltemp = #1%
   \let\@optionalnext = \relax
   \@futurenonspacelet\@optionalnext\@bracketcheck
}%
%
%
\def\@bracketcheck{%
   \ifx [\@optionalnext
      \expandafter\@@getoptionalarg
   \else
      \let\@optionalarg = \empty
      \expandafter\@optionaltemp
   \fi
}%
\def\@@getoptionalarg[#1]{%
   \def\@optionalarg{#1}%
   \@optionaltemp
}%
%
%
%
\def\@nnil{\@nil}%
\def\@fornoop#1\@@#2#3{}%
\def\@for#1:=#2\do#3{%
   \edef\@fortmp{#2}%
   \ifx\@fortmp\empty \else
      \expandafter\@forloop#2,\@nil,\@nil\@@#1{#3}%
   \fi
}%
\def\@forloop#1,#2,#3\@@#4#5{\def#4{#1}\ifx #4\@nnil \else
       #5\def#4{#2}\ifx #4\@nnil \else#5\@iforloop #3\@@#4{#5}\fi\fi
}%
\def\@iforloop#1,#2\@@#3#4{\def#3{#1}\ifx #3\@nnil
       \let\@nextwhile=\@fornoop \else
      #4\relax\let\@nextwhile=\@iforloop\fi\@nextwhile#2\@@#3{#4}%
}%
%
%
%
\@innernewif\if@fileexists
\def\@testfileexistence{\@getoptionalarg\@finishtestfileexistence}%
\def\@finishtestfileexistence#1{%
   \begingroup
      \def\extension{#1}%
      \immediate\openin0 =
         \ifx\@optionalarg\empty\jobname\else\@optionalarg\fi
         \ifx\extension\empty \else .#1\fi
         \space
      \ifeof 0
         \global\@fileexistsfalse
      \else
         \global\@fileexiststrue
      \fi
      \immediate\closein0
   \endgroup
}%
%
%
%
%
\def\bibliographystyle#1{%
   \@readauxfile
   \@writeaux{\string\bibstyle{#1}}%
}%
\let\bibstyle = \@gobble
%
%
\let\bblfilebasename = \jobname
\def\bibliography#1{%
   \@readauxfile
   \@writeaux{\string\bibdata{#1}}%
   \@testfileexistence[\bblfilebasename]{bbl}%
   \if@fileexists
      \nobreak
      \@readbblfile
   \fi
}%
\let\bibdata = \@gobble
%
%
\def\nocite#1{%
   \@readauxfile
   \@writeaux{\string\citation{#1}}%
}%
\@innernewif\if@notfirstcitation
%
%
\def\cite{\@getoptionalarg\@cite}%
%
%
\def\@cite#1{%
   \let\@citenotetext = \@optionalarg
   \printcitestart
   \nocite{#1}%
   \@notfirstcitationfalse
   \@for \@citation :=#1\do
   {%
      \expandafter\@onecitation\@citation\@@
   }%
   \ifx\empty\@citenotetext\else
      \printcitenote{\@citenotetext}%
   \fi
   \printcitefinish
}%
\def\@onecitation#1\@@{%
   \if@notfirstcitation
      \printbetweencitations
   \fi
   \expandafter \ifx \csname\@citelabel{#1}\endcsname \relax
      \if@citewarning
         \message{\@linenumber Undefined citation `#1'.}%
      \fi
      \expandafter\gdef\csname\@citelabel{#1}\endcsname{%
\strut
\vadjust{\vskip-\dp\strutbox
\vbox to 0pt{\vss\parindent0cm \leftskip=\hsize 
\advance\leftskip3mm
\advance\hsize 4cm\strut\openup-4pt 
\rightskip 0cm plus 1cm minus 0.5cm ?  #1 ?\strut}}
         {\tt
            \escapechar = -1
            \nobreak\hskip0pt
            \expandafter\string\csname#1\endcsname
            \nobreak\hskip0pt
         }%
      }%
   \fi
   \csname\@citelabel{#1}\endcsname
   \@notfirstcitationtrue
}%
%
%
\def\@citelabel#1{b@#1}%
%
%
\def\@citedef#1#2{\expandafter\gdef\csname\@citelabel{#1}\endcsname{#2}}%
%
%
%
\def\@readbblfile{%
   \ifx\@itemnum\@undefined
      \@innernewcount\@itemnum
   \fi
   \begingroup
      \def\begin##1##2{%
         \setbox0 = \hbox{\biblabelcontents{##2}}%
         \biblabelwidth = \wd0
      }%
      \def\end##1{}
      %
      %
      \@itemnum = 0
      \def\bibitem{\@getoptionalarg\@bibitem}%
      \def\@bibitem{%
         \ifx\@optionalarg\empty
            \expandafter\@numberedbibitem
         \else
            \expandafter\@alphabibitem
         \fi
      }%
      \def\@alphabibitem##1{%
         \expandafter \xdef\csname\@citelabel{##1}\endcsname {\@optionalarg}%
         \ifx\biblabelprecontents\@undefined
            \let\biblabelprecontents = \relax
         \fi
         \ifx\biblabelpostcontents\@undefined
            \let\biblabelpostcontents = \hss
         \fi
         \@finishbibitem{##1}%
      }%
      \def\@numberedbibitem##1{%
         \advance\@itemnum by 1
         \expandafter \xdef\csname\@citelabel{##1}\endcsname{\number\@itemnum}%
         \ifx\biblabelprecontents\@undefined
            \let\biblabelprecontents = \hss
         \fi
         \ifx\biblabelpostcontents\@undefined
            \let\biblabelpostcontents = \relax
         \fi
         \@finishbibitem{##1}%
      }%
      \def\@finishbibitem##1{%
         \biblabelprint{\csname\@citelabel{##1}\endcsname}%
         \@writeaux{\string\@citedef{##1}{\csname\@citelabel{##1}\endcsname}}%
         \ignorespaces
      }%
      %
      %
      \let\em = \bblem
      \let\newblock = \bblnewblock
      \let\sc = \bblsc
      \frenchspacing
      \clubpenalty = 4000 \widowpenalty = 4000
      \tolerance = 10000 \hfuzz = .5pt
      \everypar = {\hangindent = \biblabelwidth
                      \advance\hangindent by \biblabelextraspace}%
      \bblrm
      \parskip = 1.5ex plus .5ex minus .5ex
      \biblabelextraspace = .5em
      \bblhook
      \input \bblfilebasename.bbl
   \endgroup
}%
%
%
\@innernewdimen\biblabelwidth
\@innernewdimen\biblabelextraspace
%
%
%
\def\biblabelprint#1{%
   \noindent
   \hbox to \biblabelwidth{%
      \biblabelprecontents
      \biblabelcontents{#1}%
      \biblabelpostcontents
   }%
   \kern\biblabelextraspace
}%
%
%
%
\def\biblabelcontents#1{{\bblrm [#1]}}%
%
%
\def\bblrm{\rm}%
%
%
\def\bblem{\it}%
%
%
\def\bblsc{\ifx\@scfont\@undefined
              \font\@scfont = cmcsc10
           \fi
           \@scfont
}%
%
%
\def\bblnewblock{\hskip .11em plus .33em minus .07em }%
%
%
\let\bblhook = \empty
%
%
%
\def\printcitestart{[}
\def\printcitefinish{]}
\def\printbetweencitations{, }
\def\printcitenote#1{, #1}
%
%
%
\let\citation = \@gobble
%
%
%
\@innernewcount\@numparams
%
%
\def\newcommand#1{%
   \def\@commandname{#1}%
   \@getoptionalarg\@continuenewcommand
}%
%
%
\def\@continuenewcommand{%
   \@numparams = \ifx\@optionalarg\empty 0\else\@optionalarg \fi \relax
   \@newcommand
}%
%
%
\def\@newcommand#1{%
   \def\@startdef{\expandafter\edef\@commandname}%
   \ifnum\@numparams=0
      \let\@paramdef = \empty
   \else
      \ifnum\@numparams>9
         \errmessage{\the\@numparams\space is too many parameters}%
      \else
         \ifnum\@numparams<0
            \errmessage{\the\@numparams\space is too few parameters}%
         \else
            \edef\@paramdef{%
               \ifcase\@numparams
                  \empty  No arguments.
               \or ####1%
               \or ####1####2%
               \or ####1####2####3%
               \or ####1####2####3####4%
               \or ####1####2####3####4####5%
               \or ####1####2####3####4####5####6%
               \or ####1####2####3####4####5####6####7%
               \or ####1####2####3####4####5####6####7####8%
               \or ####1####2####3####4####5####6####7####8####9%
               \fi
            }%
         \fi
      \fi
   \fi
   \expandafter\@startdef\@paramdef{#1}%
}%
%
%
%
%
\def\@readauxfile{%
   \if@auxfiledone \else 
      \global\@auxfiledonetrue
      \@testfileexistence{aux}%
      \if@fileexists
         \begingroup
            \endlinechar = -1
            \catcode`@ = 11
            \input \jobname.aux
         \endgroup
      \else
         \message{\@undefinedmessage}%
         \global\@citewarningfalse
      \fi
      \immediate\openout\@auxfile = \jobname.aux
   \fi
}%
%
%
\newif\if@auxfiledone
\ifx\noauxfile\@undefined \else \@auxfiledonetrue\fi
%
%
%
%
\@innernewwrite\@auxfile
\def\@writeaux#1{\ifx\noauxfile\@undefined \write\@auxfile{#1}\fi}%
%
%
%
\ifx\@undefinedmessage\@undefined
   \def\@undefinedmessage{No .aux file; I won't give you warnings about
                          undefined citations.}%
\fi
%
%
\@innernewif\if@citewarning
\ifx\noauxfile\@undefined \@citewarningtrue\fi
%
%
%
\catcode`@ = \@oldatcatcode


\def\widestnumber#1#2{}

\def\rm{\fam0 \tenrm}

\def\fakesubhead#1\endsubhead{\bigskip\noindent{\bf#1}\par}



%
%
%

%

\font\textrsfs=rsfs10
\font\scriptrsfs=rsfs7
\font\scriptscriptrsfs=rsfs5

\newfam\rsfsfam
\textfont\rsfsfam=\textrsfs
\scriptfont\rsfsfam=\scriptrsfs
\scriptscriptfont\rsfsfam=\scriptscriptrsfs

\edef\oldcatcodeofat{\the\catcode`\@}
\catcode`\@11

\def\Cal@@#1{\noaccents@ \fam \rsfsfam #1}

\catcode`\@\oldcatcodeofat

\newpage

\head {\S1} \endhead  \resetall 
\bigskip

We continue to investigate club guessing, (see \cite[Ch.III]{Sh:g},
\cite[Ch.VI,\S2]{Sh:e}), continuing Dzamonja Shelah \cite{DjSh:691}.
\bn
\proclaim{\stag{cg.1} Claim}  Assume
\mr
\item "{$(a)$}"  $\lambda > \kappa > \theta^+$ are regular,
$\kappa \ge \sigma$
\sn
\item "{$(b)$}"  $S_\theta \subseteq \{\delta < \lambda:\text{cf}(\delta)
= \theta\}$ is stationary
\sn
\item "{$(c)$}"  $S_\kappa \subseteq \{\delta < \lambda:\text{cf}(\delta)
= \kappa\}$ is stationary
\sn
\item "{$(d)$}"  $\bar e = \langle e_\alpha:\alpha < \lambda \rangle,
e_\alpha \subseteq \alpha$ a club $(e_0 = \emptyset,e_{\alpha +1} =
\{\alpha\})$
\sn
\item "{$(e)$}"  for every club $E$ of $\lambda,(\exists^{\text{stat}} \delta
\in S_\kappa)[\delta = \sup(E \cap \text{ nacc}(e_\delta))]$
\sn
\item "{$(f)$}"  $\delta \in S_\kappa \Rightarrow (\exists^{\text{stat}} 
\alpha < \delta)(\alpha \in S_\theta \and e_\alpha = e_\delta \cap
\alpha)$
\sn
\item "{$(g)$}"  $S \subseteq \lambda$ is unbounded not reflecting
\footnote{for weaker conclusion: not reflecting outside itself, see 
\scite{gc.8}}
\sn
\item "{$(h)$}"  $\alpha \in S_\kappa \Rightarrow \text{ otp}(e_\alpha)
= \kappa$.
\ermn
\ub{Then} we can find club $c^*_\alpha$ of $\alpha$ for
$\alpha \in S_\theta$ such that
\mr
\item "{$(\alpha)$}"  for $E \subseteq \lambda$ club, 
$A \subseteq S$ unbounded we have 
$\{\delta \in S_\kappa:\{\alpha \in e_\delta \cap
S_\theta:\alpha = \sup(A \cap \text{ nacc}(c^*_\alpha))\}$ stationary in
$\delta\}$ is stationary in $\lambda$
\sn
\item "{$(\beta)$}"  $\lambda = \mu^+ \Rightarrow \text{ otp}(c^*_\alpha) <
\mu^\omega$.
\endroster
\endproclaim
\bn
Before proving
\demo{\stag{cg.2} Observation}:  Each of the statements 
$(*)_1,(*)_2,(*)_3$ implies $\otimes$ where
\mr
\item "{$\otimes$}"  there is $g:S_\theta \rightarrow \sigma$ such that
$S_{\kappa,g} \notin \text{ id}^p(\bar e \restriction S_\kappa)$ where \nl
$S_{\kappa,g} = \{\delta \in S_\kappa:\text{for every } \zeta < \sigma
\text{ we have } \{\alpha \in e_\delta:\alpha \in S_\theta \text{ and }
g(\alpha) = \zeta\}$ is stationary (in $\delta$)$\}$
\sn
\item "{$(*)_1$}"  $\sigma^+ < \kappa$
\sn
\item "{$(*)_2$}"  $\kappa^\sigma < \lambda$
\sn
\item "{$(*)_3$}"  $=(*)^3_{\lambda,\kappa,\sigma}$ there is ${\Cal P}
\subseteq [\kappa]^\sigma$ of cardinality $< \lambda$ such that for every
club $E$ of $\kappa$ for some $X \in {\Cal P}$ we have $X \subseteq E$.
\endroster
\enddemo 
\bigskip

\demo{Proof}  By guessing clubs (\cite[Ch.III,\S2]{Sh:g}) or cardinal
arithmetic, clearly 
$(*)_1 \Rightarrow (*)_3$ and $(*)_2 \Rightarrow (*)_3$ so assume
$(*)_3$.  Without loss of generality $X \in {\Cal P} \Rightarrow 
\text{ otp}(X) = \sigma$ and
let $X = \{\gamma_{X,i}:i < \sigma\},\gamma_{X,i}$ increasing with $i$;
stipulate $\gamma_{X,\sigma} = \kappa$.  For every limit ordinal $\alpha <
\lambda$ let $f_\alpha$ be an increasing continuous function from
cf$(\alpha)$ into $\alpha$ such that if $\alpha \in S_\kappa$ then 
Rang$(f_\alpha) = e_\alpha$.  For $\xi < \theta$ an $X \in {\Cal P}$ we define
$g_{X,\xi}:S_\theta \rightarrow \sigma +1$ as follows:

$$
g_{X,\xi}(\delta) = \text{ Min}\{i \le \sigma:f_{\text{otp}(e_\delta)}(\xi)
\ge \gamma_{X,i}\}
$$
\mn
and define 

$$
\align
S_{\kappa,X,\xi} = \bigl\{\delta \in S_\kappa:&\text{for every } 
i < \sigma \text{ the set}\\
  &\{\alpha \in E_\delta \cap S_\theta:g_{X,\xi}(\delta) =i\}
\text{ is stationary in } \delta \bigr\}.
\endalign
$$
\mn
Now we shall prove that:
\mr
\item "{$\boxtimes_1$}"  for some $\xi < \theta$ and $X \in {\Cal P}$ we
have $S_{\kappa,X,\xi} \notin \text{ id}^p(\bar e \restriction S_\kappa)$
i.e. for every club of $\lambda,(\exists^{\text{stat}} \delta \in
S_{\kappa,X,\xi})(\delta = \sup[E \cap \text{ nacc}(e_\delta)])$.
\ermn
Why?  If not then for every $\xi < \theta,X \in {\Cal P}$ for some club
$E_{\xi,X}$ of $\lambda,\{\delta \in S_{\kappa,X,\xi}:\delta = 
\sup(E \cap \text{ nacc}(e_\delta))\}$ is not stationary so 
without loss of generality $\delta \in S_{\kappa,X,\xi} \Rightarrow
\delta > \sup(E \cap \text{ nacc}(e_\delta))$.
Let $E =: \cap\{E_{\xi,X}:\xi < \theta \text{ and } X \in {\Cal P}\}$, so 
$E$ is a club of $\lambda$ hence for some $\delta \in S_\kappa,\delta = 
\sup(E \cap \text{ nacc}(e_\delta))$.

Let $s = \{\alpha < \kappa:\text{cf}(\alpha) = \theta \text{ and }
f_\delta(\alpha) \in S_\theta,e_{f_\delta(\alpha)} = e_\delta \cap f_\delta
(\alpha)\}$, so by an assumption, the set $s$ is stationary, so by 
Fodour and Ulam, for some $\xi < \theta$ we have:

$$
A_{\delta^*,\xi} = \{\gamma < \kappa:(\exists^{\text{stat}} \alpha \in s)
(f_\alpha(\xi) = \gamma)\}
$$
\mn
is unbounded in $\kappa$. \nl
Let $C = \{\alpha < \kappa:\alpha \text{ limit and } \alpha = \sup(\alpha
\cap A_{\delta^*,\xi}\}$, it is a club of $\kappa$.  Hence for some $X \in
{\Cal P},X \subseteq C$.  It is easy to check that $\delta$ contradicts the
choice of $E_{\xi,X}$ which is $\subseteq E$; so $\boxtimes$ really holds.
\mn
Fix $(\xi^*,X^*) \in \kappa \times {\Cal P}$ as in $\boxtimes$ and 
let $g^* = g_{\xi^*,X^*}$.  It is easy to check that $g^*$ is as required.
\hfill$\square_{\scite{cg.2}}$
\enddemo
\bigskip

\demo{Proof of \scite{cg.1}}  Let for 
$\alpha < \lambda,C_\alpha$ be a club of $\alpha$ such that
$C_{\alpha +1} = \{\alpha\},C_0 = \emptyset$ and 
$[\alpha \text{ limit } \notin A \Rightarrow
\text{ acc}(C_\alpha) \cap A = \emptyset]$.

For $\alpha \in S_\theta \cup S_\kappa$ let $e_{\alpha,n}$ be defined by

$$
e_{\alpha,0} = e_\alpha
$$

$$
\align
e_{\alpha,n+1} = e_{\alpha,n} \cup \bigl\{\beta&:\text{for some }
\gamma_0 \in e_{\alpha,n} \text{ letting } \gamma_1 = 
\text{ Min}(e_{\alpha,n} \backslash
(\gamma_0 +1)) \\
  &\text{ we have } \beta \in (\gamma_0,\gamma_1),\beta \in C_{\gamma_1}.
\bigr\}
\endalign
$$
\mn
Note
\mr
\item "{$\boxtimes_2$}"  if $\delta_1,\delta_2 \in S_\theta \cup S_\kappa,
e_{\delta_1} = e_{\delta_2} \cap \delta_1,\delta_1 \in e_{\delta_2}$ then
$\dsize \bigwedge_n[e_{\delta_1,n} = e_{\delta_2,n} \cap \delta_1]$.
\ermn
Note
\mr
\item "{$\bigodot$}"  if $A \subseteq S$ and $\beta \in A \cap \delta
\backslash e_\delta,\delta \in S_\theta$ \ub{then} for some 
$n < \omega$ we have $\beta \in \text{ nacc}(e_{\delta,n})$ \nl
[Why?  Check.]
\ermn
For $\delta \in S_\theta$ let $c^*_\delta = e_{\delta,g^*(\delta)}$ if
$g^*(\delta) < \omega$ and $e_\delta$ if $g^*(\delta) \ge \omega$.  We shall
show that $\langle c^*_\delta:\delta \in S_\theta \rangle$ is as required.  
The bound on order type if $\lambda = \mu^+$ should be clear.

Let $E$ be a club of $\lambda$ and $A \subseteq S$ be such that: $\alpha \in S
\backslash A \Rightarrow A \cap \text{ acc}(e_\alpha) = 0$;
clearly $E' = \{\delta \in E:\delta \text{ a limit ordinal } > \kappa
\text{ and otp}|A \cap \delta| = 0\}$ is a club of $E$, so
for some $\delta^* \in S_{\kappa,X^*,\xi^*}$ we have $\delta^* = \sup(E' \cap
\text{ nacc}(e_{\delta^*}))$.  By $\bigodot$ for some $n^* < \omega,\delta^*
= \sup(A \cap \text{ nacc}(e_{\delta^*,n^*}))$ hence $t = \{\alpha \in
e_{\delta^*} \cap S_\theta:e_\alpha = e_{\delta^*} \cap \alpha$ and
$g_{X^*,\xi^*}(\alpha) = n^*\}$ is stationary in $\delta^*$ hence
$t^* = \{\alpha \in t:\alpha = \sup(A \cap \text{ nacc}(e_{\delta^*,n^*}))\}$
is stationary in $\delta^*$ but for $\alpha \in t^*$ necessarily $\alpha =
\sup(A \cap \text{nacc}(e_{\delta^*,n^*})) = (A \cap \text{ nacc}
(e_{\alpha,n^*}) = \sup(A \cap \text{ nacc}(C^*_\alpha))$, so we are done.
\hfill$\square_{\scite{cg.1}}$
\enddemo
\bigskip

\proclaim{\stag{gc.4} Claim}  Let 
$\lambda,\kappa,\theta,S_\kappa,S_\theta$ satisfy clauses
(a),(b),(c),(g) from \scite{cg.1} and (f), i.e. $[\delta \in S_\kappa 
\Rightarrow \delta \cap S_\theta$ stationary in $\delta]$. \nl
\ub{Then} for some $\bar e$ the assumption (d),(e),(h) (hence the conclusion
of \scite{cg.1}) holds \ub{if} at least one of the following occurs:
\mr
\item "{$(A)$}"   $\lambda$ is a successor of regular, \nl
$\kappa < \lambda,\kappa = \text{ cf}(\kappa) > \theta + \aleph_1$,cf$(\theta)
= \theta$
\sn
\item "{$(B)$}"   $\lambda$ is a successor of regulars, $\aleph_1 = \kappa
< \lambda,\theta = \aleph_0,2^{\aleph_0} < \lambda$
\sn
\item "{$(C)$}"  $\lambda$ is a successor of a singular, $\theta < \kappa$
regular $< \lambda$ and {\rm cf}$(\lambda^-) \notin \{\theta,\kappa\}$ \nl
$S_\kappa \in I[\lambda]$ and: $\lambda$ is strong limit or at least
$(\forall \alpha < \lambda^-)(|\alpha|^\theta \le \lambda^*)$.
\endroster
\endproclaim
\bigskip

\demo{Proof}  If (A) or (B) holds, we use the square on successor of a regular
(exists by \cite[\S4]{Sh:351} or \cite[Ch.III,\S3.x]{Sh:e}) to find
$S^+ \subseteq \lambda$ on which we have square, $S^+ \cap S_\kappa$
stationary and continue as usual, (see \cite[Ch.III,2.14(3),(4)]{Sh:g} or
\cite[Ch.III,3.x]{Sh:e}) to get also club guessing.
If (C), for the existence of $\langle
e_\alpha:\alpha \in S_\theta \cup S_\kappa \rangle$ as required see (x.x)
(or just let $\lambda = \mu^+,\langle \mu_\varepsilon:\varepsilon <
\text{ cf}(\mu) \rangle$ be increasing continuous with limit $\mu,\mu_0 >
\kappa$, we can find $\bar a = \langle a_\alpha:\alpha < \lambda \rangle$
such that $a_\alpha \subseteq \alpha$, otp$(a_\alpha) \le \kappa,\beta \in
a_\alpha \Rightarrow a_\beta = a_\alpha \cap \beta$ and for some club $E$
of $\lambda$ we have $\delta \in S \cap E \Rightarrow \alpha = \sup
(a_\alpha)$). \nl
Let $C_\alpha$ be $\alpha \cap$ closure$(a_\alpha)$ for each 
$\beta > 0,\beta \in S_\theta$ let ${\Cal P}_\beta = \{c_\alpha \cap \beta:
\beta \in c_\alpha,\alpha < \lambda^+\}$ be lised as $\langle d_{\beta,\zeta}:
\zeta < \mu \rangle$, possible as $|{\Cal P}_\beta| \le \mu$. \nl
For $\alpha \in S \cap E$ let $h_\alpha:c_\alpha \cap S_\theta \rightarrow
\mu$ be defined by $h_\alpha(\beta) = \text{ Min}\{\varepsilon < \text{ cf}
(\mu):c_\alpha \cap \beta = d_{\beta,\varepsilon}\}$. \sn
Now we use the following strengthening of \cite{EK}, see ([xx])
\mr
\item "{$\otimes$}"  if cf$(\mu) \ne \theta = \text{ cf}(\theta),\lambda =
\mu^+,(\forall \alpha < \mu)(|\alpha|^\theta < \mu)$ then we can find
$f_\alpha \in {}^\lambda \mu$ for $\alpha < \mu$ such that for every $A \in
[\lambda]^\theta$ and $f \in {}^A \mu$ for some $\alpha < \mu$ we have
$|\{i \in A:f(i) = f_\alpha(i)\}| = \theta$.
\ermn
So for some $\zeta$ letting $e_\beta = \beta \cap$ closure$(d_{\beta,f_\zeta
(\beta)})$, we are done.  hfill$\square_{\scite{gc.4}}$
\enddemo
\bigskip

\proclaim{\stag{gc.5} Claim}  Assume clauses (a)-(g) of \scite{cg.1} and
\mr
\item "{$(e)^+$}"  $S_\kappa \notin \text{ id}^a(\bar e \restriction
S_\kappa)$ (easy to get, by minor corrections follows from $(e)$ by
intersecting with a suitable clue using $(h)$
\sn
\item "{$(i)$}"  $\aleph^\tau_0 \le \sigma$.
\ermn
\ub{Then} we can strengthen the conclusion in \scite{cg.2} to
\mr
\item "{$(*)_\tau$}"  if $A_\varepsilon \subseteq S$ for $\varepsilon < \tau$
and $E$ a club of $\lambda$, \ub{then} for stationarily many $\delta \in
S_\kappa$ for stationarily many $\alpha \in S_\theta \cap e_\delta$ we have
$e_\alpha = e_\delta \cap \alpha$ and $\alpha = \sup(A_\varepsilon \cap
\text{ nacc}(e^*_\alpha))$ for $\varepsilon < \tau$.
\endroster
\endproclaim
\bigskip

\demo{Proof}  The difference concerns the use of $e_{\alpha,n}$.  We define
$e_{\alpha,n}$ as before and then define $e_{\alpha,h}$ for every $\alpha \in 
S_\theta$ and $h \in {}^\tau \omega$ as follows:

$$
\align
e_{\alpha,h} = e_\alpha \cup \bigl\{ (\gamma_0,\gamma_1) \cap
e_{\alpha,n}:&\gamma_0 < \gamma_1 \text{ are successive members of} \\
  &e_\alpha \text{ and } \gamma_0 = i \text{ mod } \tau \and i < \tau
\Rightarrow n = h(i) \bigr\}.
\endalign
$$
\mn
Let $\{h_\zeta:\zeta < \aleph^\tau_0\}$ list ${}^\tau \omega$ and let
$g:S_\theta \rightarrow \sigma$ be as in the previous proof and lastly
$c^*_\alpha = e_{\alpha,(h_{g(\alpha)})}$.  The rest should be clear.
\hfill$\square_{\scite{gc.5}}$
\enddemo
\bigskip

\remark{\stag{gc.7} Remark}  1) If we restrict ourselves to 
stationary $A \subseteq S$,
then for club $E$ we can demand $\alpha = \sup(A \cap E \cap \text{ nacc}
(c^*_\alpha))$ as we can work on $A$. 
\endremark
\bigskip

\proclaim{\stag{gc.8} Claim}  Assume
\mr
\item "{$(a)$}"  $\lambda > \kappa > \theta$ are regular, $\kappa > \sigma
\ge \sigma_1 \ge \aleph^\tau_0$
\sn
\item "{$(b)$}"  $S_\theta \subseteq \{\delta < \lambda:\text{cf}(\delta)
= \theta\}$ stationary
\sn
\item "{$(c)$}"  $S_\kappa \subseteq \{\delta < \lambda:\text{cf}(\delta) =
\kappa$ and $S_\theta \cap \delta$ is stationary$\}$ is stationary
\sn
\item "{$(d)$}"  $\bar e^* = \langle e^*_\alpha:\alpha < \lambda \rangle,
e^*_\alpha \subseteq \alpha,e^*_{\alpha +1} = \{\alpha\},\alpha$ limit
$\Rightarrow e^*_\alpha$ a club of $\alpha$
\sn
\item "{$(e)$}"  $S \subseteq \lambda$ and $\bar e^*$ runs away from $S$ which
means $\alpha \in \lambda \backslash S \and \alpha$ limit $\Rightarrow S \cap
\text{ acc}(e^*_\alpha) = \emptyset$ so this implies: 
$S$ does not reflect outside itself
\sn
\item "{$(f)$}"  $(*)^3_{\lambda,\sigma,\sigma_1}$ from \scite{cg.2}.
\ermn
\ub{Then} we can find $\bar c^* = \langle c^*_\alpha:\alpha \in S_\theta
\rangle$ such that
\mr
\item "{$(\alpha)$}"  $c^*_\alpha$ a club of $\alpha$
\sn
\item "{$(\beta)$}"  if $A_\varepsilon \subseteq S$ is unbounded for
$\varepsilon < \tau$ and runs away from $\bar e^*$, \ub{then} for
stationarily many $\delta \in S_\kappa$, for stationarily many $\alpha <
\delta$, we have $\dsize \bigwedge_{\varepsilon < \tau} \alpha = \sup
(A_\varepsilon \cap \text{ nacc}(c^*_\alpha))$.
\endroster
\endproclaim
\bigskip

\demo{Proof}  For limit $\alpha < \lambda$ let $\langle \gamma_{\alpha,
\varepsilon}:\varepsilon < \text{ cf}(\alpha) \rangle$ be increasing
continuous with limit $\alpha$.

By Fodor and Ulam for every $\delta \in S_\kappa$ there are $\xi_\delta <
\theta$ and $\gamma_\delta < \delta$ such that
\mr
\item "{$\bigodot_1$}"  cf$(\gamma_\delta) = \sigma$
\sn
\item "{$\bigodot_2$}"  for every $\gamma < \gamma_\delta$ the set \nl
$\{\alpha \in S_\theta \cap \delta:\gamma_{\alpha,\xi} \in (\gamma,
\gamma_\delta)\}$ is stationary in $\delta$.
\ermn
By Fodor for some $\gamma^*,\xi^*$ we have

$$
S^*_\kappa = \{\delta \in S_\kappa:\xi_\delta = \xi^* \text{ and }
\gamma_\delta = \gamma^*\} \text{ is stationary in } \lambda.
$$
\mn
For every $X \in {\Cal P}$ (from $(*)^3_{\lambda,\sigma,\sigma_1}$) letting
$X = \{\beta_{X,i}:i < \sigma_1\}$, we let 
$g_X:S_\theta \rightarrow \sigma_1+1$
be $g_X(\alpha) = i$ iff $\gamma_{\alpha,\xi^*} \in [\beta_{X,i},\beta
_{X,i+1})$ and $\sigma_1$ if no such $X$. \nl
Now for some $X = X^*$
\mr
\item "{$\bigodot_3$}"  for every club $E$ of $\lambda$ for stationarily
many $\delta \in S^*_\kappa$, for every $i < \sigma_1$ for stationarily many
$\alpha \in S_\theta \cap \delta$ (stationary in $\delta$!) we have
$g_{X^*}(\alpha)=i$.
\ermn
Now we continue as before.  \hfill$\square_{\scite{gc.8}}$
\enddemo

\bigskip\bigskip\bigskip
    
REFERENCES.  
\bibliographystyle{lit-plain}
\bibliography{lista,listb,listx,listf,liste}

\enddocument

\bye